\documentclass{amsart}
\usepackage{amscd,amssymb,subfigure}
\usepackage[arrow,matrix,graph]{xy}

\newtheorem{theorem}{Theorem}[section]
\newtheorem{corollary}[theorem]{Corollary}
\newtheorem{lemma}[theorem]{Lemma}
\newtheorem{prop}[theorem]{Proposition}

\theoremstyle{definition}
\newtheorem{definition}[theorem]{Definition}

\newtheorem{remark}[theorem]{Remark}
\newtheorem*{ack}{Acknowledgment}
\newtheorem*{add}{Added in Proof}

\newtheorem{claim}{Claim}

\numberwithin{equation}{section}

\newcommand{\A}{{\mathcal A}}
\newcommand{\B}{{\mathcal B}}
\renewcommand{\L}{{\mathcal L}}
\newcommand{\X}{{\mathcal X}}
\newcommand{\oA}{\A^*}

\newcommand{\LL}{{\mathbf L}}

\newcommand{\Z}{\mathbb{Z}}
\newcommand{\Q}{\mathbb{Q}}

\newcommand{\C}{\mathbb{C}}

\renewcommand{\k}{\Bbbk}
\newcommand{\HH}{{\mathfrak H}}

\newcommand{\G}{\mathsf{G}}
\newcommand{\E}{\mathsf{E}}
\newcommand{\V}{\mathsf{V}}
\newcommand{\T}{\mathsf{T}}

\DeclareMathOperator{\rank}{rank}
\DeclareMathOperator{\gr}{gr}

\DeclareMathOperator{\id}{id}

\DeclareMathOperator{\ch}{char}
\DeclareMathOperator{\Lie}{{Lie}}

\newcommand{\surj}{\twoheadrightarrow}
\newcommand{\abs}[1]{\left| #1 \right|}

\begin{document}

\title[Decomposition of Lie algebras associated to arrangements]{%
When does the associated graded Lie algebra of an
arrangement group decompose?}

\author[Stefan Papadima]{Stefan Papadima$^1$}
\address{Institute of Mathematics of the Academy,
P.O. Box 1-764,
RO-014700 Bucharest, Romania}
\email{Stefan.Papadima@imar.ro}

\author[Alexander~I.~Suciu]{Alexander~I.~Suciu$^2$}
\address{Department of Mathematics,
Northeastern University,
Boston, MA 02115, USA}
\email{a.suciu@neu.edu}
\urladdr{http://www.math.neu.edu/\~{}suciu}

\thanks{$^1$Partially supported by CERES grant 152/2003 of
the Romanian Ministry of Education and Research}

\thanks{$^2$Partially supported by NSF grant DMS-0311142}



\begin{abstract}
Let $\A$ be a complex hyperplane arrangement, with
fundamental group $G$ and holonomy Lie algebra $\HH$.
Suppose  $\HH_3$ is a free abelian group of
minimum possible rank, given the values the M\"obius
function $\mu\colon \L_2\to \Z$ takes on the rank $2$ flats of $\A$.
Then the associated graded Lie algebra of $G$ decomposes
(in degrees $\ge 2$) as a direct product of free Lie algebras.
In particular, the ranks of the lower central series quotients of
the group are given by
$\phi_r(G)=\sum _{X\in \L_2} \phi_r(F_{\mu(X)})$, for $r\ge 2$.
We illustrate this new Lower Central Series formula with
several families of examples.
\end{abstract}
\maketitle

\section{Introduction}
\label{sect:intro}

\subsection{}
The purpose of this paper is to give an answer to the question
posed in the title.  Let $\A$ be an arrangement of finitely
many hyperplanes through the origin of $\C^{\ell}$, and denote by
$G(\A)=\pi_1\big(\C^{\ell}\setminus \bigcup_{H\in \A} H\big)$ the
fundamental group of its complement. In Section \ref{sect:state},
we single out a class of arrangements, closely related to certain 
arrangements studied in \cite{CS}, \cite{SS}.  Roughly speaking, 
$\A$ is  {\em decomposable} if a certain quadratic,
graded Lie algebra $\HH(\A)$, naturally defined in terms of the
codimension $2$ flats of $\A$, has minimal possible dimension
in degree $3$, over any ground field.

Our main result (Theorem~\ref{thm:intro}) implies the following:  
If $\A$ is decomposable, then the associated graded Lie algebra 
of $G(\A)$ decomposes as a direct product of free Lie algebras 
(in degrees $r\ge 2$):
\begin{equation}
\label{grA}
\gr_{\ge 2}(G(\A)) \cong \prod_{X\in \L_2(\A)} \gr_{\ge 2}(F_{\mu(X)}).
\end{equation}
Here:
\begin{itemize}
\item $\L(\A)=\left\{X=\bigcap_{H \in \B}H \mid \B\subseteq \A\right\}$
is the intersection lattice, $\L_2(\A)$
is the set of codimension~$2$ flats, and $\mu\colon \L(\A)\to \Z$
is the M\"obius function.
\item $\{\Gamma_r G\}_{r\ge 1}$ is the lower
central series, given by $\Gamma_1 G=G$ and
$\Gamma_{r+1} G=(\Gamma_{r} G,G)$.
\item $\gr(G)=\bigoplus_{r\ge 1}  \Gamma_{r} G /\Gamma_{r+1} G$,
with Lie bracket induced by the group  commutator.
\item $F_n$ is the free group of rank $n$, and
$\gr(F_n)=\LL_n$ is the free Lie algebra on $n$ generators.
\end{itemize}
Moreover, as we show in Proposition~\ref{prop:heredity}, the 
decomposability property of $\A$ is inherited by all 
sub-arrangements of $\A$.

\subsection{}
The associated graded Lie algebra $\gr(G(\A))$ is not {\it a priori}\/
determined by the intersection lattice, and, as such, it is not easy 
to handle.  We turn instead to a more manageable, combinatorial 
approximation:  The {\em holonomy Lie algebra} of the arrangement, 
$\HH(\A)$, defined as the quotient of $\LL(\A)$, the free Lie algebra 
on variables $\{x_H \mid H\in \A\}$, modulo the ideal $J(\A)$ generated 
by relations corresponding to rank $2$ flats:
\begin{equation}
\label{holoA}
\HH(\A)=\LL( \A)\Big\slash \text{ideal}\, \Big\{
\big[x_H , \sum_{H'\in \A \colon H'\supset X} x_{H'}\big] \:\big| \:
X \in \L_2(\A) \text{ and } X\subset H \big. \Big\}  \Big. .
\end{equation}

As shown by Kohno \cite{Ko} (based on foundational work by
Sullivan \cite{Su} and Morgan \cite{Mo}), the associated graded
Lie algebra $\gr(G(\A))$ and the holonomy Lie algebra $\HH(\A)$
are rationally isomorphic:
\begin{equation}
\label{kohno}
\gr(G(\A))\otimes \Q\cong \HH(\A)\otimes \Q.
\end{equation}
At the integral level, there is a surjective
Lie algebra map, $\Psi_{\A}\colon \HH(\A)\surj \gr(G(\A))$, 
such that $\Psi_{\A}\otimes \Q$ is an isomorphism, 
see \cite{MP}.  In general, there exist arrangements 
for which $\Psi_{\A}$ is not injective.  Nevertheless, 
for the class of decomposable arrangements we consider
here, $\Psi_{\A}$ gives an isomorphism $\gr(G(\A))\cong \HH(\A)$,
see Theorem \ref{thm:intro}\eqref{t2}.

\subsection{}
The lower central series ranks of a finitely-generated group 
$G$ are defined as $\phi_r(G)=\rank \gr_r(G)$.
For a free group, the LCS ranks are given by
Witt's formula:  $\prod_{r=1}^{\infty}(1-t^r)^{\phi_r(F_n)}=1-n t$.
For an arrangement group, the LCS ranks are determined
by the intersection lattice, via \eqref{kohno} and \eqref{holoA}.
Clearly, $\phi_1(G(\A))=\abs{\A}$.  Hence, to determine the 
LCS ranks of $G(\A)$, we only need to compute the graded 
ranks of the derived holonomy algebra,
$\HH'(\A)=\bigoplus_{r\ge 2} \HH_r(\A)$.

From work of Falk \cite{Fa}, we know that 
$\dim_{\Q} \HH_r(\A)\otimes \Q \ge 
\sum_{X\in \L_2(\A)} \phi_r(F_{\mu(X)})$, for all $r\ge 2$,
with equality holding for $r=2$. Guided by these facts, we say that
$\A$ is decomposable if the lower bound is attained in degree $3$, 
for every field $\k$:
\begin{equation}
\label{eq=decintro}
\dim_{\k} \HH_3(\A) \otimes \k = 
\sum_{X\in \L_2(\A)} \phi_3(F_{\mu(X)}).
\end{equation}

Up to now, an explicit formula for the LCS ranks of an 
arrangement group has only been known in the case 
when the  intersection lattice is supersolvable \cite{FR}, 
or, more generally, hypersolvable \cite{JP}.  The 
isomorphism \eqref{grA} leads to a new LCS formula, for 
the combinatorially defined class of decomposable arrangements:
\begin{equation}
\label{newlcs}
\prod_{r=1}^{\infty}(1-t^r)^{\phi_r(G(\A))}=
(1-t)^{\abs{\A}}\prod_{X\in \L_2(\A)} 
\frac{1- \mu(X) t}{(1-t)^{\mu(X)}}\, .
\end{equation}

This LCS formula verifies the more general ``resonance 
LCS formula," conjectured in \cite{S}, in what is arguably 
the simplest, yet most basic case.

As a byproduct of our main theorem, we compute in 
Section~\ref{sect:chen} the integral {\em Chen Lie algebra}\/ 
of a decomposable arrangement, and we also obtain the 
Chen analog of decomposition \eqref{grA}, thus improving 
upon results from \cite{CS}.

\subsection{}
Formula \eqref{newlcs} is equivalent to
$\phi_r(G(\A))=\sum _{X\in \L_2(\A)} \phi_r(F_{\mu(X)})$,
for all $r\ge 2$. In other words, the (higher) LCS ranks 
behave as if $G(\A)$ were to decompose as a direct product 
of free groups, of ranks dictated by the M\"obius function.
This happens, for instance, for the class of (hypersolvable, 
decomposable) arrangements considered in \cite{CDP}, 
where the arrangement group is always a product of free 
groups. In general, though, the group of a decomposable 
arrangement does not decompose in this manner.

\begin{figure}[ht]
\subfigure{%
\begin{minipage}[t]{0.35\textwidth}
\setlength{\unitlength}{0.45cm}
\begin{picture}(3,3.3)(-2,-0.2)
\put(3,3){\line(1,-1){3}}
\put(3,3){\line(-1,-1){3}}
\put(3,3){\line(0,-1){3}}
\put(1.5,1.5){\line(1,0){3}}
\put(0,0){\line(1,0){6}}
\multiput(0,0)(6,0){2}{\circle*{0.5}}
\multiput(1.5,1.5)(1.5,0){3}{\circle*{0.5}}
\multiput(3,3)(0,-3){2}{\circle*{0.5}}
\end{picture}
\end{minipage}
}
\subfigure{%
\begin{minipage}[t]{0.35\textwidth}
\setlength{\unitlength}{0.45cm}
\begin{picture}(3,3.3)(-2.8,-0.2)
\put(3,3){\line(1,-1){3}}
\put(3,3){\line(-1,-1){3}}
\put(0,0){\line(1,0){6}}
\multiput(0,0)(6,0){2}{\circle*{0.5}}
\multiput(1.5,1.5)(3,0){2}{\circle*{0.5}}
\multiput(3,3)(0,-3){2}{\circle*{0.5}}
\end{picture}
\end{minipage}
}
\caption{\textsf{The $\operatorname{X}_2$ and 
$\operatorname{X}_3$ matroids}}
\label{fig:x3x2}
\end{figure}

For example, consider  the $\operatorname{X}_2$ and 
$\operatorname{X}_3$ arrangements, whose matroids 
are depicted in Figure~\ref{fig:x3x2}.  It is readily checked that both 
arrangements are decomposable (compare with \cite{CS}, \cite{SS}), 
but not hypersolvable (see Remark~\ref{rem=grlike}).  

For the $\operatorname{X}_2$ arrangement, we find that 
$\phi_r(G(\A))= \phi_r((F_2)^{\times 5})$, for all $r\ge 2$,
yet $\phi_1(G(\A))<\phi_1((F_2)^{\times 5})$; thus, 
$G(\A)\not\cong (F_2)^{\times 5}$.

For the $\operatorname{X}_3$ arrangement, we find that   
$\phi_r(G(\A))= \phi_r((F_2)^{\times 3})$, for all $r\ge 1$.  
Even so, $G(\A)\not\cong (F_2)^{\times 3}$.  Indeed, it 
can be checked that $G(\A)\cong G\times \Z$, where $G$ is 
the celebrated Stallings group,\footnote{As a consequence, 
we can compute the LCS ranks of the Stallings group: 
$\phi_1(G)=5$, and $\phi_r(G)= \phi_r((F_2)^{\times 3})$, 
for  $r\ge 2$.  In \cite{PS-bb}, we give an LCS formula 
that applies to {\em any}\/ Bestvina-Brady group 
associated to a connected flag complex.} 
equal to the kernel of the projection 
$(F_2)^{\times 3} \to \Z$, which sends each 
standard generator to $1$.  As shown in \cite{St}, 
the group $H_3(G)$ is not finitely generated.  
It follows that $G(\A)$ does not admit a finite 
$K(G(\A),1)$; in particular, $G(\A)$ cannot be 
isomorphic to {\em any}\/ finite direct product of  
free groups of finite rank.

In view of these examples, and of the infinite families 
of decomposable, non-hyper\-solvable graphic arrangements 
from Section~\ref{sect:graphic}, we see that the LCS 
formula \eqref{newlcs} is a genuinely new formula,
with a range of applicability which overlaps only marginally 
with that of the classical LCS  formula.

\begin{ack}
Most of this work was done while the second author visited the
Institute of Mathematics of the Romanian Academy in June--July 2003,
with partial support from grant CERES/CE4 of the Romanian Ministry of
Education and Research.
\end{ack}

\section{Decomposable arrangements}
\label{sect:state}

In this section, we delineate the class of decomposable arrangements,
and state in detail our main result.

\subsection{}
\label{subsect:proj}
For an arrangement $\A$, denote by $\Z^{\A}$ the
free abelian group on $\A$, with basis $\{x_H\mid H\in \A\}$.
For a sub-arrangement $\B\subset \A$, 
let $\pi_{\B}\colon \Z^{\A} \to \Z^{\B}$
be the canonical projection map, defined by $\pi_{\B}(x_H)=x_H$,
if $H\in \B$, and $\pi_{\B}(x_H)=0$, if $H\notin \B$, and
let $\LL (\pi_{\B})\colon \LL(\A)\to \LL(\B)$ be its extension
to free Lie algebras.  Clearly, $\LL (\pi_{\B})(J(\A))\subset J(\B)$,
and so we get a Lie algebra epimorphism,
\begin{equation}
\label{hhpi}
\HH (\pi_{\B})\colon \HH(\A) \surj \HH(\B).
\end{equation}

For a flat $X\in \L_2(\A)$, let $\A_X=\{ H\in \A \mid H\supset X\}$
be the {\em localization} of $\A$ at $X$.
This is a pencil of $\abs{\A_X}=\mu(X)+1$ hyperplanes.
The group $G(\A_X)$ is isomorphic to $F_{\mu(X)}\times \Z$;
thus, $\gr(G(\A_X)) \cong \LL_{\mu(X)} \times \LL_1$. 
From the defining relations \eqref{holoA}, we also have
$\HH(\A_X) \cong \LL_{\mu(X)} \times \LL_1$, and so
$\HH(\A_X)\cong\gr(G(\A_X))$.

Set $\pi_X=\pi_{\A_X}$.
The maps $\HH(\pi_X)\colon \HH(\A)\to \HH(\A_X)$ assemble
into a Lie algebra map from $\HH(\A)$ to the direct product
of the holonomy Lie algebras of its localized sub-arrangements:
\begin{equation}
\label{pimap}
\pi=\left(\HH(\pi_X)\right)_X\colon \HH(\A)\longrightarrow
\prod_{X\in \L_2(\A)} \HH(\A_X).
\end{equation}

The starting point of our investigation is the following result,
to be proved in \S\ref{subsect:surjproof}.

\begin{prop}
\label{prop:surj}
The restriction of $\pi$ to derived subalgebras,
\[
\pi'\colon \HH'(\A)\rightarrow \prod_{X \in \L_2(\A)} \HH'(\A_X), 
\]
is surjective.
\end{prop}

By comparing ranks of the source and target of
$\pi_r\colon \HH_r(\A)\rightarrow \prod_{X} \HH_r(\A_X)$ for $r\ge 2$, 
we recover a lower bound for the LCS ranks of an arrangement group,
first obtained by M.~Falk \cite{Fa}, by other methods.

\begin{corollary}[\cite{Fa}]
\label{cor:phibound}
For all $r\ge 2$,
\begin{equation}
\label{phibound}
\phi_r(G(\A))\ge \sum_{X\in \L_2(\A)} \phi_r(F_{\mu(X)}).
\end{equation}
\end{corollary}

\subsection{}
\label{subsect:decomp}
Our main goal here is to understand when the natural map $\pi'$ from
Proposition \ref{prop:surj} is, in fact, an isomorphism; in particular,
when the inequalities \eqref{phibound} become equalities.

It is easy to see that $\pi_2$ is always an isomorphism. On the
other hand, the maps $\pi_r$ ($r\ge 3$) may not be isomorphisms,
as illustrated by the braid arrangements $\B_{\ell}$ in $\C^{\ell}$
($\ell\ge 4$).  In this case, the LCS formula of Kohno \cite{Ko85}
and Falk-Randell \cite{FR}, when applied to the pure braid group
$P_{\ell}=G(\B_{\ell})$, shows that inequality \eqref{phibound} is
strict in degree $r= 3$.

This prompts the following definition.

\begin{definition}
\label{def:kdecomp}
Let $r\ge 2$ be an integer, and let $\k$ be a field.
We say that $\HH_r(\A)$ is {\em $\k$-decomposable} if
\begin{equation}
\label{eq:kdecomp}
\dim_{\k} \HH_r(\A) \otimes \k = \sum_{X\in \L_2(\A)} \phi_r(F_{\mu(X)}).
\end{equation}
An arrangement $\A$ is {\em decomposable} if $\HH_3(\A)$ is 
$\k$-decomposable, for every field $\k$.
\end{definition}

By Proposition \ref{prop:surj}, $\HH_r(\A)$ is $\k$-decomposable
if and only if $\pi_r\otimes \k$ is an isomorphism, whereas $\A$
is decomposable precisely when $\pi_3$ is an isomorphism.

\subsection{}
\label{subsect:h34}
Two other decomposability conditions were considered in \cite{CS} 
and \cite{SS}.  Let us briefly compare those conditions to ours.

The condition from \cite{CS} entails the decomposability of the
$I$-adic completion of the Alexander invariant of $G(\A)$ as the direct
sum of the $I$-adic completions of the Alexander invariants of
$G(\A_X)$, taken over $X\in \L_2(\A)$.  It can be shown that this
condition on Alexander invariants is equivalent, over $\Q$, 
to the decomposability of $\HH_3(\A)$, in the sense of
Definition~\ref{def:kdecomp}.

The condition from \cite{SS} entails the minimality of the linear
strand of the free resolution of the Orlik-Solomon algebra of $\A$ as a
module over the corresponding exterior algebra. As stated in
\cite[Definition 2.10]{SS}, the MLS condition is equivalent to the
$\k$-decomposability of $\HH_3(\A)$, for $\k$ a field of characteristic $0$.
Actually, the only place where the hypothesis $\ch\k=0$ is needed
in that context is to insure that $\dim_{\k}\gr_*(G(\A))\otimes \k =
\dim_{\k}\HH_*(\A)\otimes \k$. All the other homological algebra
arguments work as well over a field of positive characteristic.
Consequently, Theorem~5.6 from \cite{SS} gives the following:
If $\HH_3(\A)$ is $\k$-decomposable, then $\HH_4(\A)$ is 
$\k$-decomposable.  In particular, if $\A$ is decomposable 
(i.e., $\pi_3$ is an isomorphism), then $\pi_4$ is an isomorphism.

\subsection{}
\label{subsect:main}
Our main result is Theorem \ref{thm:intro} below, which improves
upon the aforementioned  result from \cite{SS}, in several ways.
For one, it pushes the range where $\pi_r$ is an isomorphism
from $r=4$ to infinity.  For another, it assembles the graded pieces
$\pi_r$ ($r\ge 2$) into a Lie algebra isomorphism between the derived
holonomy Lie algebra of $\A$ and a product of derived free Lie algebras.
Finally, it gives a new LCS-type formula for the group of a decomposable
arrangement, thus verifying Conjecture 5.7 from \cite{SS}.

\begin{theorem}
\label{thm:intro}
Let $\A$ be a decomposable arrangement. Then:
\begin{enumerate}
\item \label{t1}
$\gr(G(\A))\cong\HH(\A)$, as graded Lie algebras.
\item \label{t2}
$\HH(\A)$ is torsion-free, as a graded abelian group.
\item  \label{t3}
$\pi'\colon \HH'(\A)\longrightarrow \prod_{X\in \L_2(\A)} \HH'(\A_X)$
is an isomorphism of graded Lie algebras.
\item  \label{t4} The LCS ranks $\phi_r=\phi_r(G(\A))$ are given
by the following combinatorial formula:
\begin{equation}
\label{declcs}
\prod_{r=1}^{\infty}(1-t^r)^{\phi_r}=
(1-t)^{b_1-b_2}\prod_{X\in \L_2(\A)} (1- \mu(X) t),
\end{equation}
where $b_1=\abs{\A}$ and $b_2=\sum_{X\in \L_2(\A)} \mu(X)$.

\end{enumerate}
\end{theorem}

Here is an immediate corollary, already mentioned in the Introduction.
\begin{corollary}
\label{cor=z}
If $\A$ is decomposable,
then the associated graded Lie algebra of $G(\A)$ decomposes
as a direct product of free Lie algebras (in degrees $r \ge 2$):
\begin{equation}
\label{grdec}
\gr_{\ge 2}(G(\A)) \cong \prod_{X\in \L_2(\A)} \gr_{\ge 2}(F_{\mu(X)}).
\end{equation}
\end{corollary}

Over the rationals, we can be even more precise:  
$\HH_3(\A)$ is decomposable over $\Q$, i.e.,
$\phi_3(G(\A))=\sum _{X\in \L_2(\A)} \phi_3(F_{\mu(X)})$,
if and only if the derived subalgebra of the rational 
associated graded Lie algebra of $G(\A)$ decomposes 
as a direct product of derived free Lie algebras over $\Q$:
\begin{equation}
\label{qgrdec}
\left(\gr(G(\A))\otimes \Q\right)' \cong 
\prod_{X\in \L_2(\A)} \LL_{\mu(X)}'\otimes \Q.
\end{equation}
This follows from Proposition~\ref{iotaprime} below and 
formula \eqref{kohno}.

Now set $\phi_r^{\k}(\A):= \dim_{\k} \HH_r(\A)\otimes \k$, 
for $\k$ a field and $r\ge 1$. The isomorphism \eqref{kohno} 
implies $\phi_r^{\Q}(\A)= \phi_r(G(\A))$, for all $r$. 
As a consequence of 
Theorem \ref{thm:intro}, we obtain the following 
characterization of decomposability, in terms of 
LCS-type formulas in arbitrary characteristic.

\begin{corollary}
\label{cor=lcschar}
The arrangement $\A$ is decomposable if and only if,
for every field $\k$, 
\[
\prod_{r=1}^{\infty}(1-t^r)^{\phi^{\k}_r(\A)}=
(1-t)^{b_1-b_2}\prod_{X\in \L_2(\A)} (1- \mu(X) t).
\]
\end{corollary}

Finally, suppose $\A$ is hypersolvable, with exponents 
$ d_1=1, d_2, \dots, d_{\ell}$.  The Poincar\'{e} polynomial 
of the quadratic Orlik-Solomon algebra associated to $\A$ 
is then given by $\overline{P}_{\A}(t)=\prod_{i=1}^{\ell}(1+ d_i t)$; 
see \cite[Proposition 3.2]{JP}.  
Putting together the decomposable LCS formula \eqref{declcs} and
the hypersolvable LCS formula from \cite[Theorem C]{JP}, we obtain
the following relationship between the exponents $d_i$ and the 
level-$2$ M\"obius function $\mu\colon \L_2(\A)\to \Z$ of a 
decomposable, hypersolvable arrangement $\A$. 
(We will exploit this relationship in the last section, 
within the framework of graphic arrangements.)

\begin{corollary}
\label{cor=dechs}
If $\A$ is both hypersolvable and decomposable, then:
\[
\prod_{i=1}^{\ell}(1+ d_i t)=
(1+t)^{\abs{\A}}\prod_{X\in \L_2(\A)} \frac{1+ \mu(X) t}{(1+t)^{\mu(X)}}.
\]
\end{corollary}

\section{The $\iota$ map}
\label{sect:iota}

In this section, we define the natural candidate for the inverse map to
$\pi'\colon \HH'(\A) \to \prod_{X\in \L_2(\A)} \HH'(\A_X)$, and discuss
some of its properties.

\subsection{} 
\label{subsect:iotab}

Let $\B$ be a sub-arrangement of $\A$.
Let $\iota_{\B}\colon \Z^{\B} \to \Z^{\A}$ be the canonical
inclusion, defined by $\iota_{\B}(x_H)=x_H$, and let
$\LL (\iota_{\B})\colon \LL(\B)\to \LL(\A)$ be its extension to free
Lie algebras.  In general, the map  $\LL (\iota_{\B})$ need not preserve
the defining ideals of the holonomy Lie algebras of $\B$ and $\A$.

However, suppose $\B$ is {\em closed} in $\A$, i.e., the only linear
combinations of defining forms for the hyperplanes in $\B$ which
are defining forms for hyperplanes in $\A$ are (up to constants)
the defining forms for the hyperplanes in $\B$.
Then $\L_2(\B)=\{ X \in \L_2(\A) \mid \A_X \subset \B\}$.
Thus, $\LL (\iota_{\B})(J(\B))\subset J(\A)$, and so we
get a map of graded Lie algebras,
\begin{equation}
\label{hhiota}
\HH (\iota_{\B})\colon \HH(\B) \to \HH(\A).
\end{equation}

For a flat $X\in \L_2(\A)$, note that $\A_X$ is closed in $\A$.
Set $\iota_X=\iota_{\A_X}$.

\begin{lemma}
\label{lem:ortho}
Let $X, Y\in \L_2(\A)$, and $\B\subset \A$.  Then:
\begin{enumerate}
\item \label{i1}
$\HH(\pi_X)\circ \HH(\iota_X) = \id$.
\item \label{i2}
$\HH'(\pi_{\B})\circ \HH'(\iota_X) =0$, if $\abs{\B\cap \A_X}\le 1$.
\item \label{i3}
$\HH'(\pi_X)\circ\HH'(\iota_Y) = 0$ if $X\ne Y$.
\end{enumerate}
\end{lemma}

\begin{proof}
\eqref{i1}
Clearly, $\pi_X\circ \iota_X$ is the identity map on $\Z^{\A_X}$.

\eqref{i2}
For each $r\ge 2$, the group $\HH_r(\A_X)$ is generated by elements of 
the form $x=[x_{H_1},[x_{H_2}, \cdots [x_{H_{r-1}},x_{H_r}]\cdots]]$, 
where $H_1,\dots, H_r$ are hyperplanes in $\A_X$.  Now, since 
$\abs{\B\cap \A_X}\le 1$, one of those hyperplanes, say $H_i$, must 
not belong to $\B$; otherwise, $H_1= H_2= \dots =H_r$, and so $x=0$.  
Hence, by definition, $\pi_{\B}(x_{H_i})=0$, and so 
\[
\HH(\pi_{\B})\circ \HH(\iota_X ) (x) =
[\pi_{\B}(x_{H_1}),[ \pi_{\B}(x_{H_2}), \cdots
[\pi_{\B}(x_{H_{r-1}}),\pi_{\B}(x_{H_r})]\cdots ]] = 0.
\]
Thus, $\HH(\pi_{\B})\circ \HH(\iota_X )=0$ in degrees $\ge 2$.

\eqref{i3}
If $X\ne Y$, then $\abs{\A_X\cap \A_Y} \le 1$.
Hence \eqref{i2} applies.
\end{proof}

\subsection{Proof of Proposition~\ref{prop:surj}}
\label{subsect:surjproof}

The maps $\HH(\iota_X)$ define
a homomorphism of graded abelian groups,
\begin{equation}
\label{iota}
\iota\colon  \prod_{X\in \L_2(\A)} \HH(\A_X)\longrightarrow \HH(\A).
\end{equation}
Let $\iota'\colon  \prod_{X} \HH'(\A_X)\rightarrow \HH'(\A)$ be the
restriction of $\iota$ to derived subalgebras. The orthogonality relations from
Lemma~\ref{lem:ortho} imply that $\pi'\circ \iota'=\id$.
Thus, $\pi'$ is surjective, and so Proposition~\ref{prop:surj} 
is proved. \hfill\qed

\subsection{}  
\label{subsect:linalg}

The following Lemma (the proof of which is an exercise in linear algebra)
will be used repeatedly later on.

\begin{lemma}
\label{linalg}
Let $U$ and $\{V_X\}_{X\in \X}$ be finite-dimensional vector spaces over
a field $\k$, with $\abs{\X}$ finite.  Set  $V=\bigoplus_{X} V_X$.
Suppose we have linear maps
$\pi_X\colon U\to V_X$ and $\iota_X\colon V_X \to U$ such that
$\pi_X\circ \iota_Y=\delta_{X,Y}$. Set $\pi=(\pi_X)_{X}\colon U\to V$
and $\iota=\sum_X \iota_X\colon V \to U$.  Then, the
following conditions are equivalent:
\begin{enumerate}
\item\label{i} $\pi$ is an isomorphism.
\item\label{ii} $\iota$ is surjective.
\item\label{iii} $\sum_X \iota_X\circ \pi_X = \id_U$.
\item\label{iv} $\dim_{\k} U = \sum_X \dim_{\k} V_X$.
\end{enumerate}
\end{lemma}

\subsection{} 
\label{subsect:heredity}

By Lemmas \ref{lem:ortho} and \ref{linalg}, $\HH_3(\A)$ is
$\k$-decomposable if and only if the map
\begin{equation}
\label{iota3}
\iota_3\otimes \k\colon \bigoplus_{X\in \L_2(\A)}
\HH_3(\A_X)\otimes \k \longrightarrow \HH_3(\A) \otimes \k
\end{equation}
is surjective.  We use this criterion to show that decomposability is hereditary.

\begin{prop}
\label{prop:heredity}
If $\B$ is a sub-arrangement of $\A$, and
if $\HH_3(\A)$ is $\k$-decomposable, then $\HH_3(\B)$ is also
$\k$-decomposable.
\end{prop}

\begin{proof}
Note that $\L_2(\B)=\left\{ X  \in \L_2(\A) \mid
\abs{\A_X \cap \B } \ge 2 \right\}$. Furthermore, if $X\in \L_2(\B)$,
then $\B_X=\A_X\cap \B$.
Consider the following diagram:
\begin{equation}
\label{subi}
\xymatrix{%
\bigoplus_{X\in \L_2(\A)} \HH'(\A_X) \ar[rr]^(0.65){\iota'^{\A}}
 \ar@{>>}[d]^{\rho}&& \HH'(\A) \ar@{>>}[d]^{\HH'(\pi^{\A}_{\B})}\\
\bigoplus_{X\in \L_2(\B)} \HH'(\B_X) \ar[rr]^(0.65){\iota'^{\B}} &&
\HH'(\B)
}
\end{equation}
where $\rho$ restricts to
$\HH'\big(\pi^{\A_X}_{\B_X}\big)\colon \HH'(\A_X) \to \HH'(\B_X)$
if $X\in \L_2(\B)$, and $\rho=0$ otherwise.
Diagram \eqref{subi} commutes. Indeed, if $X\in \L_2(\B)$, this is clear.
If $X\notin\L_2(\B)$, then $\abs{\A_X \cap \B} \le 1$,
and so, by Lemma~\ref{lem:ortho}\eqref{i2},
$\HH'\big(\pi^{\A}_{\B}\big)\circ \HH'\big(\iota^{\A}_{\A_X}\big)=0$.

Now, if $\HH_3(\A)$ is $\k$-decomposable, then $\iota_3^{\A}\otimes \k$
is surjective. From the commutativity of diagram \eqref{subi}, we infer that
$\iota_3^{\B}\otimes \k$ is also surjective, and we are done.
\end{proof}

\section{Surjectivity of $\iota'$}
\label{sect:iotaprime}

In general, the map
$\iota'\colon  \prod_{X} \HH'(\A_X)\rightarrow \HH'(\A)$
is not surjective.  On the other hand, if $\A$ is decomposable,
$\iota'$ is surjective.  This we show in the next Proposition, which is
the key to our main result.

\begin{prop}
\label{iotaprime}
Suppose $\HH_3(\A)$ is $\k$-decomposable. Then
\[
\iota_r \otimes \k\colon  \prod_X \HH_r (\A_X)\otimes \k
\longrightarrow \HH_r (\A)\otimes \k
\]
is surjective, for all $r\ge 2$.
\end{prop}

\begin{proof}
For simplicity, we will suppress the field $\k$ from the notation.
We will need to establish various commutation relations in $\HH(\A)$, between
elements in $\HH(\iota_X) (\HH(\A_X))$ and $\HH(\iota_Y) (\HH(\A_Y))$,
where  $X$ and $Y$ are distinct flats in $\L_2(\A)$.  Again for simplicity,
we will suppress the inclusion $\iota$ from the notation, and work in $\HH(\A)$.

Since $X\ne Y$, there are two possibilities:  either $\A_X \cap \A_Y=\emptyset$,
or $\A_X\cap \A_Y$ consists of a single hyperplane.  Pick $H'\in \A_X$
and $H''\in \A_Y$, so that, if $\oA_X=\A_X \setminus \{H'\}$ and
$\oA_Y=\A_Y \setminus \{H''\}$ are the corresponding deletions, then
\begin{equation}
\label{cap}
\oA_X \cap \A_Y = \A_X \cap \oA_Y =\emptyset.
\end{equation}

Let us note the following fact, whose proof is immediate, and which will be
used repeatedly in the sequel.
For any flat  $Z\in \L_2(\A)$, and for any hyperplane $H\in \A_Z$,
\begin{equation}
\label{lie}
\HH'(\A_Z)=\Lie^{>1}(\oA_Z),
\end{equation}
where $\oA_Z=\A_Z\setminus \{H\}$, and where
$\Lie^r(\oA_Z)$ denotes the degree $r$ piece of
the Lie sub\-algebra generated by $\{x_K \mid K\in \oA_Z\}$
inside $\HH(\A_Z)$.

Here is the first commutation property, for which the decomposability
assumption on $\HH_3(\A)$ is needed in a crucial way.

\begin{claim}
\label{step1}
If $H_i\in \oA_X$ and $c\in \HH_2(\A_Y)$, then $[x_{H_i},c]=0$.
\end{claim}

\begin{proof}
By \eqref{lie}, it is enough to verify the claim for $c\in  \Lie^2(\oA_Y)$.
Apply $\HH(\pi_{Z})$, for some $Z\in \L_2(\A)$.  If $Z\ne Y$, we get
$[\HH(\pi_{Z})(x_{H_i}),\HH(\pi_Z)(c)]=[\HH(\pi_{Z})(x_{H_i}),0]=0$.
If $Z =Y$, we get $[\HH(\pi_{Y})(x_{H_i}),c]=[0,c]=0$, since
$\oA_X \cap \A_Y = \emptyset$. Thus, $\pi_3([x_{H_i},c])=0$,
and so $[x_{H_i},c]=0$, since, by assumption, $\pi_3$ is an isomorphism.
\end{proof}

Using Claim~\ref{step1}, we obtain the next commutation property.

\begin{claim}
\label{step2}
If $b\in \HH_2(\A_X)$ and $c\in \HH_s(\A_Y)$ ($s\ge 2$), then $[b,c]=0$.
\end{claim}

\begin{proof}
As before, we may assume that $b\in \Lie^2 (\oA_X)$ and $c\in \Lie^s (\oA_Y)$.
The proof is by induction on $s$. For $s=2$, Claim~\ref{step2}
follows from Claim~\ref{step1}, via the Jacobi identity.
 For the induction step, take an element
$c\in \Lie^{s+1}(\oA_Y)$,
and write it as $c=[x_{H_j},c']$, with $H_j \in \oA_Y$ and $c'\in \Lie^s(\oA_Y)$.
By the Jacobi identity,
\[
[b,c] = [[b,x_{H_j}],c']+[x_{H_j},[b,c']].
\]
Note that $[b,x_{H_j}]=0$ by Claim~\ref{step1}, and $[b,c']=0$ by induction.
Thus, $[b,c]=0$.
\end{proof}

Finally, using both Claims~\ref{step1} and~\ref{step2}, we prove the following
key commutation property.

\begin{claim}
\label{step3}
If $H_i\in \oA_X$ and $c\in \HH_{s}(\A_Y)$ ($s\ge 2$), then
$[x_{H_i},c]\in \HH_{s+1}(\A_Y)$.
\end{claim}

\begin{proof}
The proof is by induction on $s$.  The case $s=2$ follows from
Claim~\ref{step1}.  For the induction step, take an element
$c=[x_{H_j},c']\in \HH_{s+1}(\A_Y)$, with
$H_j\in \A_Y$ and $c'\in \HH_s(\A_Y)$. By the Jacobi identity,
\[
[x_{H_i},c] = [[x_{H_i},x_{H_j}],c']+[x_{H_j},[x_{H_i},c']].
\]
By induction, $[x_{H_i},c']\in \HH_{s+1}(\A_Y)$, and so
$[x_{H_j},[x_{H_i},c']]\in \HH_{s+2}(\A_Y)$.

On the other hand, since $H_i\ne H_j$,
there is a flat $Z\in \L_2(\A)$ such that $\{H_i,H_j\}\subset \A_Z$,
and so $[x_{H_i},x_{H_j}]\in \HH_2(\A_Z)$.
If $Z=Y$, then $[[x_{H_i},x_{H_j}],c']\in
[\HH_2(\A_Y),\HH_s(\A_Y)]\subset \HH_{s+2}(\A_Y)$.
If $Z\ne Y$, then $[[x_{H_i},x_{H_j}],c']=0$, by Claim~\ref{step2}.
Either way, we conclude that $[x_{H_i},c]\in \HH_{s+2}(\A_Y)$.
\end{proof}

Having established the above claims,
we are now ready to prove Proposition~\ref{iotaprime},  by induction on $r$.
For $r=2$, the map $\iota_2$ is surjective, since $\pi_2\circ \iota_2=\id$,
and $\pi_2$ is an isomorphism (for arbitrary $\A$).
For the induction step, it is plainly enough to show that
\begin{equation}
\label{eq=ionto}
[x_H, c]\in \HH_{r+1}(\A_Y),
\end{equation}
for any $H\in \A$ and $c\in \HH_r(\A_Y)$, where $Y\in \L_2(\A)$ and $r\ge 2$.

If $H\in \A_Y$, this is clear. Assuming $H\notin \A_Y$, pick any
$H''\in \A_Y$, and set $X=H\cap H''\in \L_2(\A)$. Note that $X\neq Y$,
and $\oA_X =\A_X \setminus \{H''\}$. Hence, $H\in \oA_X$, and \eqref{eq=ionto}
now follows from Claim~\ref{step3}. The proof of Proposition~\ref{iotaprime}
is thus complete.
\end{proof}

\section{Proof of Theorem~\ref{thm:intro}}
\label{proofthm}
We are now in position to prove our main result.

Let $\A$ be an arbitrary arrangement.  Recall we defined in
\S\ref{subsect:proj} a homomorphism of graded Lie algebras,
$\pi\colon \HH(\A)\to \prod_{X\in \L_2(\A)} \HH(\A_X)$.
Recall also we defined in \S\ref{subsect:surjproof}  a homomorphism of graded
abelian groups, $\iota\colon  \prod_{X\in \L_2(\A)} \HH(\A_X)\to \HH(\A)$,
with the property that $\pi'\circ \iota'=\id$, which showed that
$\pi'$ is an epimorphism.

Now suppose $\A$ is decomposable.
By Proposition \ref{iotaprime}, each map $\iota_r$ ($r\ge 2$)
is surjective.  By Lemmas~\ref{lem:ortho} and~\ref{linalg}, each map
$\pi_r$  ($r\ge 2$) is an isomorphism.
Hence, $\pi'$ is an isomorphism of Lie algebras,
with inverse $\iota'$.  This proves Part \eqref{t3} of Theorem~\ref{thm:intro}.

Part \eqref{t2} follows at once from \eqref{t3}, and the fact that each Lie
algebra $\HH(\A_X)\cong \LL_{\mu(X)}\times \LL_1$ is torsion-free.

Part \eqref{t1} follows from \eqref{t2}, together with \eqref{kohno}
and \cite[Proposition~5.1]{MP}.

Part \eqref{t4}  follows from \eqref{t3}, together with \eqref{kohno} and
the discussion from \cite[\S1.5]{SS}. \hfill\qed

\section{Decomposable Chen Lie algebras}
\label{sect:chen}

Another, much coarser approximation to the associated graded
Lie algebra of a group is its Chen Lie algebra.
We now study the effect of the decomposability
condition on the Chen Lie algebra of an arrangement group.

Given a finitely-generated group $G$, let $G/G''$ be the quotient by
its second derived subgroup.  We call the associated graded Lie algebra
$\gr(G/G'')$, the {\em Chen Lie algebra} of $G$.
Set  $\theta_k(G)=\rank \gr_k(G/G'')$.  Plainly, $\theta_k(G)=\phi_k(G)$
for $k\le 3$, and $\theta_k(G)\le \phi_k(G)$ for $k> 3$.

Now suppose $G(\A)$ is an arrangement group.  Then, as shown
in \cite[Theorem 11.1]{PS04}, there is an isomorphism of graded Lie algebras,
\begin{equation}
\label{eq=chenf}
\gr(G(\A)/G''(\A))\otimes \Q \cong (\HH(\A)/\HH''(\A)) \otimes \Q.
\end{equation}
Let $B(\A)=\HH'(\A)/\HH''(\A)$ be the {\em infinitesimal Alexander invariant}
of $\A$.
Taking graded ranks on both sides of \eqref{eq=chenf}, we find:
\begin{equation}
\label{thb}
\theta_k(G(\A)) = \rank B_k(\A), \quad
\text{for all $k\ge 2$}.
\end{equation}

Recall once more the surjective map of graded Lie algebras
from Proposition~\ref{prop:surj},
$\pi'\colon \HH'(\A)\surj \prod_{X\in \L_2(\A)} \HH'(\A_X)$.
By abelianization, we obtain an
epimorphism of graded abelian groups,
\begin{equation}
\label{bp}
B(\pi)\colon B(\A) \longrightarrow
\bigoplus_{X\in \L_2(\A)} B(\A_X)\,.
\end{equation}
By comparing graded ranks of the source and target of $B(\pi)$,
we recover a lower bound for the Chen ranks of an arrangement group,
first obtained in \cite{CS} by other methods.

\begin{corollary}[\cite{CS}]
\label{cor:thetabound}
For all $r\ge 2$,
\begin{equation}
\label{thetabound}
\theta_r(G(\A))\ge \sum_{X\in \L_2(\A)} \theta_r(F_{\mu(X)}),
\end{equation}
where $\theta_r(F_n)=(r-1) \binom{n+r-2}{r}$.
\end{corollary}

Note that $B_2(\pi)=\pi_2$ is an isomorphism, and thus equality holds in
\eqref{thetabound} for $r=2$.  For $r\ge 3$, though, the inequality can
well be strict; see again \cite{CS}.

As another application of our methods, we provide a complete
description of the Chen Lie algebra of a decomposable arrangement.

\begin{theorem}
\label{tmm:chendeco}
If $\A$ is decomposable, then:
\begin{enumerate}
\item $\gr(G(\A)/G''(\A))=\HH(\A)/\HH''(\A)$,
as graded Lie algebras over $\Z$.
\item $\gr(G(\A)/G''(\A))$ is torsion-free, as a graded abelian group.
\item The Chen ranks of $G(\A)$, for $r\ge 2$, are given by
\begin{equation}
\label{thetadec}
\theta_r(G(\A))= \sum_{X\in \L_2(\A)} \theta_r(F_{\mu(X)}).
\end{equation}
\end{enumerate}
\end{theorem}

\begin{proof}
For any flat $X\in \L_2(\A)$, we have
$B(\A_X)=\LL'_{\mu(X)} / \LL''_{\mu(X)}$,
which is known to be torsion-free.
Now, since $\A$ is decomposable,
Theorem~\ref{thm:intro}\eqref{t3} implies that
$B(\pi)$ is an isomorphism, and consequently
$B(\A)$ is torsion-free, as well.
Hence, $\HH(\A)/\HH''(\A)$ is also torsion-free.
Parts (1) and (2) now follow from Theorem B in \cite{PS04}.
Part (3) follows from the fact that $B(\pi)$ is an isomorphism,
and \eqref{thb}.
\end{proof}

Formula \eqref{thetadec} was derived by other methods in \cite{CS},
under the decomposability condition from that paper.

\section{Decomposable graphic arrangements}
\label{sect:graphic}

To a (simple) graph $\G$, with vertex set $\V=\{1,\dots,\ell\}$ and
edge set $\E$, there corresponds a {\em graphic arrangement}
in $\C^{\ell}$, denoted by $\A_{\G}$.
The hyperplane corresponding to an edge $e=(i,j)$
is $H_e=\{z_i-z_j=0\}$.
For example, if $\G=K_{\ell}$, the complete graph on $\ell$ vertices,
then $\A_{K_{\ell}}=\B_{\ell}$, the braid arrangement in $\C^{\ell}$.

For each flat $X\in \L_2(\A_{\G})$, there are either $2$ or $3$
hyperplanes containing $X$. Under the identification $\A_{\G}=\E$,
a flat of size $3$ of corresponds to a triangle in the graph,
while a flat of size $2$ corresponds to a pair of edges which
is not included in any element of the triangle-set $\T$.
Thus, the holonomy Lie algebra of $\A_{\G}$ can be identified with
the quotient of the free Lie algebra on variables $e\in \E$ by
the corresponding ideal of quadratic relations:
\begin{equation}
\label{holograph}
\HH(\G)=\LL(\E)\Big\slash \text{ideal} \left\{
\begin{array}{ll}
\left[e_1 , e_2 + e_3 \right], & \text{if\,  $\{e_1,e_2,e_3\}\in \T$}
\\
\left[ e_1, e_2 \right] , &\text{if\, $\{e_1,e_2,e\} \notin \T,\:
\forall e\in \E$}
\end{array}
\right\} \Big. .
\end{equation}

As shown in \cite{SS}, the $\Q$-decomposability condition for a graphic
arrangement can be read off the graph itself, as the absence
of complete quadrangles in $\G$.  We present a strengthened form
of this result, which nicely illustrates our methods.

\begin{prop}
\label{prop:graphdec}
For a graphic arrangement $\A_{\G}$, the following conditions are
equivalent:
\begin{enumerate}
\item \label{gr1}  $\A_{\G}$ is decomposable.
\item \label{gr2}  $\HH_3(\G)$ is decomposable over some field $\k$.
\item \label{gr3}  $\G$ contains no complete subgraphs on $4$ vertices.
\end{enumerate}
\end{prop}

\begin{proof}
The implication \eqref{gr1} $\Rightarrow$ \eqref{gr2} is obvious.

To show \eqref{gr2} $\Rightarrow$ \eqref{gr3}, suppose $K_4$
is a subgraph of $\G$.  Then, the braid arrangement $\B=\A_{K_4}$
is a sub-arrangement of $\A_{\G}$. But $\B$
is not $\k$-decomposable, for any field $\k$.
Indeed, $\sum_{X\in \L_2(\B)} \phi_3(F_{\mu(X)})=8$,
whereas $\dim_{\k} \HH_3(\B)\otimes \k=10$
(see \cite{Ko85} for the case $\k=\Q$, and \cite{JP} for the
general case).
This contradicts Proposition~\ref{prop:heredity}.

To show \eqref{gr3} $\Rightarrow$ \eqref{gr1}, we must check
that $\HH_3(\G)$ is spanned by $\{ \iota_{\tau} (\HH_3(\tau)) \mid \tau\in \T\}$,
where $\iota_{\tau}\colon \HH(\tau) \to \HH(\G)$ is the natural
inclusion.  As an abelian group,  $\HH_3(\G)$ is generated by
elements of the form $x=[e_1,[e_2,e_3]]$.
Note that the edges $e_2,e_3$ must belong to a common triangle,
say, $\tau$, for, otherwise, $[e_2,e_3]=0$ in $\HH(\G)$.
If  $e_1 \in \tau$, then clearly $x\in \iota_{\tau}(\HH_3(\tau))$.
If $e_1\notin \tau$, we will show that $x=0$, and that will finish
the proof.

First, we claim that there are two edges, $e$ and $e'$,
in $\tau$ such that
\begin{equation}
\label{ee}
[e_1,e]=[e_1,e']=0.
\end{equation}
To verify the claim, denote by $\G_0$ the subgraph supported on the
vertices of $\tau$ and $e_1$.
Since $e_1 \notin \tau$, there are two possibilities:
\begin{enumerate}
\item[(a)] $\G_0$ has $5$ vertices.  Then $[e_1, e]=0$, for any edge
$e\in \tau$.
\item[(b)] $\G_0$ has $4$ vertices.
Since by assumption $\G_0\ne K_4$, again there are two possibilities:
\begin{enumerate}
\item[(b$_1$)]  $\G_0$ has $4$  edges.  
Then $[e_1,e]=0$, for any edge $e\in \tau$.
\item[(b$_2$)] $\G_0$ has $5$  edges. 
Then $\G_0$ is the union of two triangles,
with an edge in common.  If $e, e'$ are the other 
two edges in $\tau$, then $[e_1,e]=[e_1,e']=0$.
\end{enumerate}
\end{enumerate}
Thus, \eqref{ee} holds in all cases.

Now, applying \eqref{lie} to $\tau^*=\{e,e'\}$, we see that
$[e_2,e_3]$ is a multiple of $[e,e']$ in $\HH_2(\tau)$.
Hence, by \eqref{ee} and the Jacobi identity, $x=0$.
\end{proof}

Let $\kappa_s=\kappa_s(\G)$ be the number of $K_{s+1}$
subgraphs of $\G$; for example, $\kappa_0=\abs{\V}$,
$\kappa_1=\abs{\E}$, $\kappa_2=\abs{\T}$.  If $\kappa_3=0$,
then, by Proposition~\ref{prop:graphdec} and Theorem~\ref{thm:intro},
we have
\begin{equation}
\label{graphlcs}
\prod_{r=1}^{\infty} (1-t^r)^{\phi_r} =
(1-t)^{\kappa_1-2\kappa_2} (1-2t)^{\kappa_2}.
\end{equation}

We now provide concrete examples where this LCS formula
gives new information. For that, we need graphs which are
not chordal (i.e., supersolvable), or, more generally,
hypersolvable (in the sense of \cite{PS02}), since otherwise, 
previously known formulas apply.

\begin{prop}
\label{prop:nonhs}
Let $\G$ be a graph with $\kappa_1\le 2\kappa_2$ and $\kappa_3=0$.
Then $\A_{\G}$ is decomposable, but not hypersolvable.
\end{prop}

\begin{proof}
Since $\kappa_3=0$, the arrangement $\A_{\G}$ is decomposable.
If $\A_{\G}$ were hypersolvable, then, by the LCS formula from \cite{JP},
$\prod_{r=1}^{\infty} (1-t^r)^{\phi_r} = (1-t) P(t)$, for some
polynomial $P$. In view of \eqref{graphlcs}, this can only happen
when $\kappa_1-2\kappa_2>0$.
\end{proof}

\begin{remark}
\label{rem=grlike}
Let $\A$ be an arrangement (not necessarily graphic) for which the
M\"obius function takes only the values $1$ and $2$ on $\L_2(\A)$.
Set
\[
\kappa_1 =\abs{\A},\quad \text{and} \quad
\kappa_2 =\abs{\{ X\in \L_2(\A) \mid
\mu (X)=2 \}} .
\]
The same argument as in Proposition~\ref{prop:nonhs}
shows the following:  If $\A$ is decomposable and hypersolvable, 
then $\kappa_1 - 2\kappa_2 >0$.
\end{remark}

\begin{figure}[ht]
\subfigure{%
\begin{minipage}[t]{0.35\textwidth}
\setlength{\unitlength}{0.35cm}
\begin{picture}(3,3.7)(-3,-1.7)
\xygraph{!{0;<8mm,0mm>:<0mm,8mm>::}
[]*-{\bullet}
(-[dd]*-{\bullet}
,-[dr]*R(10){\G}*L(12.5){\leadsto}*-{\bullet}
(-[ur],-[dl]),-[rr]*-{\bullet}
(-^{e}[dd]*-{\bullet}
(-[ll],-[ul])))
}
\end{picture}
\end{minipage}
}
\subfigure{%
\begin{minipage}[t]{0.35\textwidth}
\setlength{\unitlength}{0.35cm}
\begin{picture}(3,3.7)(-2.7,-1.7)
\xygraph{!{0;<8mm,0mm>:<0mm,8mm>::}
[]*-{\bullet}
(-[dd]*-{\bullet},-[dr]*R(7.5){\G'}*-{\bullet}
(-[ur],-[dl]),
-[rr]*-{\bullet}
(-^{f_1}[dr]*L(3){\scriptstyle{w}}*-{\bullet}
,-^{e}[dd]*-{\bullet}
(-[ll],-[ul],-_{f_2}[ur])))
}
\end{picture}
\end{minipage}
}
\caption{\textsf{Coning an edge}}
\label{fig:graph}
\end{figure}

Now suppose $\G$ is a graph with $\kappa_1- 2\kappa_2\le 0$ 
and $\kappa_3=0$.  One can create a new graph, $\G'$, with 
the same properties, as follows.   Choose an edge $e$ of $\G$, 
pick a new vertex $w$, and join it by edges $f_1$ and $f_2$ 
to the endpoints of $e$, as in Figure~\ref{fig:graph}.
Clearly, $\V'=\V\cup \{w\}$,
$\E'=\E\cup \{f_1,f_2\}$, $\T'=\T\cup\{e,f_1,f_2\}$, and 
there are no complete quadrangles introduced.  Thus,
$\kappa_1'-2\kappa_2'=(\kappa_1+2)-2(\kappa_2+1)\le 0$, 
and $\kappa_3'=0$.  Moreover, it is easy to check that 
$\A_{\G}$ is solvable in $\A_{\G'}$, in the sense of \cite{JP}.

This permits us to create infinite families of graphs satisfying
the hypothesis of Proposition~\ref{prop:nonhs}.
For instance, start with the graph $\G^0=\G$ from the
above figure (see also \cite[Example~6.14]{SS}), and
define inductively a sequence of graphs $\{\G^i\}$ by
$\G^{i}=(\G^{i-1})'$.
Since $\G$ satisfies $\kappa_1-2\kappa_2=\kappa_3=0$,
all the graphic arrangements $\A_{\G^{i}}$ are
decomposable, but not hypersolvable.  By
\eqref{graphlcs}, the LCS ranks of the corresponding
arrangement groups are given by
$\prod_{r=1}^{\infty} (1-t^r)^{\phi_r} = (1-2t)^{i+4}$.

\begin{remark}
To the best of our knowledge, the decomposable arrangements
discussed in this paper provide the first non-hypersolvable 
examples where the LCS ranks $\phi_r$ are computed for 
{\em all}\/ values of $r$.  Note that the two graphic arrangements 
in Examples 3.7 and 5.4 from \cite{Pe} are hypersolvable.   
Indeed, the two underlying graphs can be obtained from 
hypersolvable graphs, by iterating the above construction: 
the first one, starting from a $4$-cycle, and the second 
one, starting from the $K_4$ graph.  As such, both 
arrangements are hypersolvable, cf.~\cite[\S 6]{PS02}.
The first one has rank $4$ and exponents $\{1,1,1,1,2\}$;  
the second one has rank $6$ and exponents $\{1,2,2,2,2,3\}$, 
and thus is actually supersolvable (by \cite[Theorem D]{JP}), 
despite a claim to the contrary in \cite{Pe}.
\end{remark}

\begin{add}
Using the holonomy Lie algebra approach, 
P.~Lima-Filho and H.~Schenck have recently announced 
in \cite{LS}  a proof of the LCS formula for graphic 
arrangements, as conjectured in \cite{SS}. 
\end{add}

\end{document}